\documentclass[final,leqno,onefignum,onetabnum]{siamltexmm}

\usepackage{amsmath}
\usepackage{amsfonts}

\usepackage{subcaption}

\usepackage{epsfig}

\title{Data-Driven Reduction for Multiscale Stochastic Dynamical Systems \thanks{
C.J.D. was supported by the Department of Energy Computational Science Graduate Fellowship (CSGF), grant number DE-FG02-97ER25308, and the National Science Foundation Graduate Research Fellowship, Grant No. DGE 1148900.
R.T. was supported by the European Union's Seventh Framework Programme (FP7) under Marie Curie Grant 630657 and by the Horev Fellowship.
R.T and R.R.C. were supported by the National Science Foundation, Award No. 1309858.
I.G.K. was supported by the NSF and the AFOSR.}}

\author{Carmeline~J.~Dsilva\footnotemark[2] \and Ronen~Talmon\footnotemark[3] \and C.~William~Gear\footnotemark[2] \and Ronald~R.~Coifman\footnotemark[4] \and Ioannis~G.~Kevrekidis\footnotemark[2]\ \footnotemark[5]}

\begin{document}
\maketitle
\newcommand{\slugmaster}{%
\slugger{siads}{xxxx}{xx}{x}{x--x}}

\renewcommand{\thefootnote}{\fnsymbol{footnote}}

\footnotetext[2]{Department of Chemical and Biological Engineering, Princeton University, Princeton, New Jersey, 08544, USA}
\footnotetext[3]{Department of Electrical Engineering, Technion – Israel Institute of Technology, Haifa, Israel, 3200003}
\footnotetext[4]{Department of Mathematics, Yale University, New Haven, Connecticut, 06520, USA}
\footnotetext[5]{Program in Applied and Computational Mathematics, Princeton University, Princeton, New Jersey, 08544, USA}

\renewcommand{\thefootnote}{\arabic{footnote}}

\begin{abstract}

Multiple time scale stochastic dynamical systems are ubiquitous in science and engineering,
and the reduction of such systems and their models to only their slow components is often
essential for scientific computation and further analysis.
Rather than being available in the form of an explicit analytical model, often
such systems can only be observed as a data set which exhibits dynamics on several time scales.
We will focus on applying and adapting data mining and manifold learning techniques to detect the slow components in such multiscale data.
Traditional data mining methods are based on metrics (and thus, geometries) which are
not informed of the multiscale nature of the underlying system dynamics;
such methods cannot successfully recover the slow variables.
Here, we present an approach which utilizes  both the local geometry and the {\em local dynamics} within the data set through a
metric which is both insensitive to the fast variables and more general than simple statistical averaging.
Our analysis of the approach provides conditions for successfully recovering the underlying slow variables, as well as an empirical protocol guiding the selection of the method parameters.

\end{abstract}



\pagestyle{myheadings}
\thispagestyle{plain}
\markboth{C.~J. DSILVA {\it ET AL}}{DATA-DRIVEN REDUCTION OF SDES}

\section{Introduction}

Dynamical systems of engineering interest often contain several disparate time scales.
When the evolving variables are strongly coupled, resolving the
dynamics at all relevant scales can be computationally challenging and pose problems for analysis.
Often, the goal is to write a reduced system of equations which accurately captures the dynamics on the slow time scales.
These reduced models can greatly accelerate simulation efforts, and are more appropriate for integration into larger modeling frameworks.

Following the methods of Mori \cite{mori1965transport}, Zwanzig \cite{zwanzig1961memory}, and others \cite{brey1981nonlinear, chorin2000optimal, hijon2010mori}, one can reduce the number of variables needed to describe a system of differential equations.
However, in general, this reduction introduces memory terms.
It transforms a system of differential equations into a system of (lower-dimensional) {\em integro-}differential equations,
so that the reduction of the number of variables is counterpoised by the increased complexity of the  reduced model.
Here, we will study the special case of evolution equations which contain an inherent time scale separation;
in this case, it is possible, in principle, to obtain a reduced system of differential equations in only the slow variables {\em without memory terms}.
Such an analysis crucially hinges on knowing in which variables (or, more generally, functions of variables) one can write such a reduced system of slow evolution equations.

\begin{figure}[t]

\centering

\begin{subfigure}{0.4\textwidth}
\centering
\epsfig{width=2in, file=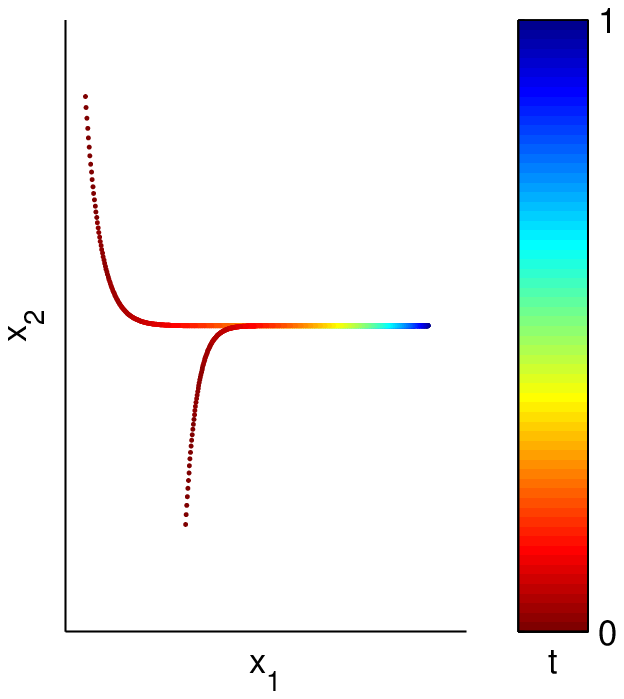}
\caption{}
\end{subfigure}
\begin{subfigure}{0.4\textwidth}
\centering
\epsfig{width=2in, file=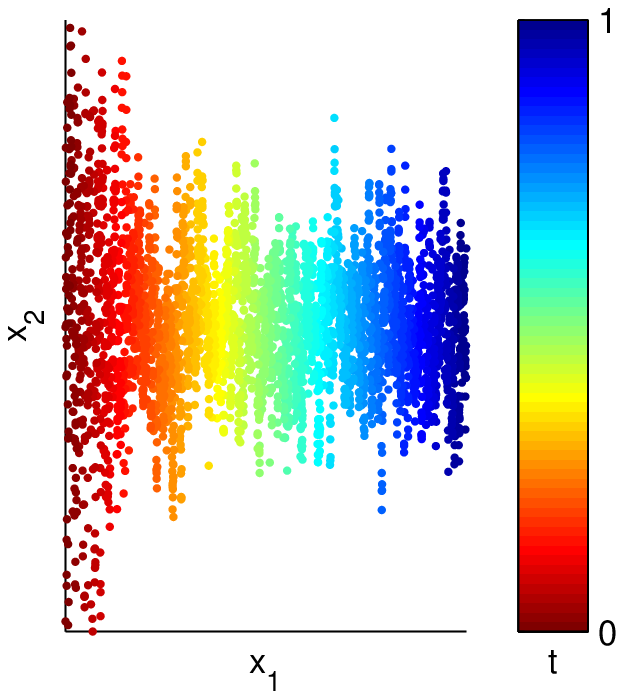}
\caption{}
\end{subfigure}
\caption{(a) Schematic of a two-dimensional, two-time scale ($\tau_1 = 50$ and $\tau_2=2$) ordinary differential equation system where the value of $x_2$ becomes slaved to the value of $x_1$.
In such an example, traditional data mining algorithms are sufficient to recover the slow variable. (b) Schematic of a two-scale two-dimensional stochastic dynamical system where {\em the statistics} of $x_2$ become slaved to $x_1$.
In such an example, traditional data mining algorithms will not recover the slow variable if the variance in the fast variable is too large. }
\label{fig:schematic_fastslow}
\end{figure}

Moving averages and subsampling have often been used in simple cases as appropriate functions of variables in which to formulate slow lower-dimensional models \cite{pavliotis2007parameter}.
However, if the underlying dynamics are sufficiently nonlinear, such statistics may
fail to capture the relevant structures and time scales within the data (see Figure~\ref{fig:schematic_fastslow} for a schematic illustration).
For well-studied systems, one often has some {\em a priori} knowledge of the appropriate observables (such as phase field variables) with which to formulate the reduced dynamics \cite{chen2002phase, wheeler1992phase}.
However, such observables may not be immediately obvious upon inspection for new complex systems, and so we require an automated approach to construct such slow variables.

Given an explicit system of ordinary differential equations,
one can make numerical approximations, such as the quasi-steady state approximation \cite{segel1989quasi} or
the partial equilibrium approximation \cite{gallagher1986combined}, to reduce the system dimensionality without introducing memory terms.
There has been some recent analytical work on extending and
generalizing such ideas to more complex systems of equations \cite{ait2008closed, calderon2007fitting, contou2011model, dong2007simplification, givon2004extracting, pavliotis2007parameter,  sotiropoulos2009model}.
However, in many instances, closed form, analytical models are not given explicitly,
but can only be inferred from simulation and/or experimental data.
We therefore turn to data-driven techniques to analyze such systems and uncover the relevant dynamical modes.
In particular, we will use a manifold-learning based approach, as such methods can accommodate nonlinear structures in high-dimensional data.

The core of most manifold learning methods is having a notion of similarity between data points,
usually through a distance metric \cite{Belkin2003, Coifman2006, coifman2005geometric, roweis2000nonlinear, tenenbaum2000global}.
The distances are then integrated into a global parametrization of the data, typically through the solution of an eigenproblem.
In this paper, we will analyze multiple time scale stochastic dynamical systems using data-driven methods.
Standard ``off-the-shelf'' manifold learning techniques which utilize the Euclidean distance are not appropriate
for analyzing data from such multiscale systems, since this metric does not account for the disparate time scales.
Research efforts have addressed the construction of more informative distance metrics, which are less sensitive to noise
and can better recover the true underlying structure in the data by suppressing unimportant sources of variability \cite{berry2013time, gepshtein2013image, rubner2000earth, simonyan2013fisher, xing2002distance}.
The Mahalanobis distance is one such metric.
It was shown that the Mahalanobis distance can remove the effect of {\em observing}
the underlying system variables through a complex, nonlinear function \cite{dsilva2013nonlinear, singer2008non, talmon2013empirical}.
Here, we will show the analogy between removing the effects of such nonlinear observation functions (in the context of data analysis), and reducing a dynamical system to remove the effects of the fast variables.
Our approach will build a parametrization of the data which is consistent with the underlying slow variables.
Because our approach is data-driven, we require no explicit description of the model, and can extract the underlying slow variables
from either simulation or experimental data.
Furthermore, the approach implicitly identifies the slow variables within the data and
does not require any {\em a priori} knowledge of the fast or slow variability sources.
Even when the underlying dynamical system is complex with nonlinear coupling between the fast and slow variables, we will show that our approach has the potential to isolate the underlying slow modes.

We will present detailed analysis for our method, and provide conditions under which it will successfully recover the slow variables.
Furthermore, based on this analysis, we will present data-driven protocols to tune the parameters of the method appropriately.
Our presentation and discussion will address two-time-scale stochastic systems; however, we claim that
our framework and analysis readily extends to systems with multiple time scale separations.

\section{Multiscale Stochastic Systems} \label{subsec:multiscale_SDE}

Consider the following two-time-scale system of stochastic differential equations (SDEs),
\begin{equation} \label{eq:general_SDE}
\begin{aligned}
dx_i(t) &= a_i(\mathbf{x}(t)) dt + dW_i(t), & \: 1 \le i \le m \\
dx_i(t) &= \frac{a_i(\mathbf{x}(t))}{\epsilon} dt + \frac{1}{\sqrt{\epsilon}} dW_i(t) , & \: m+1 \le i \le n
\end{aligned}
\end{equation}
where $W_i(t)$ are independent standard Brownian motions, $\mathbf{x}(t)  = \begin{bmatrix} x_1(t) & \cdots & x_n(t) \end{bmatrix}^T \in \mathbb{R}^n$, and $\epsilon \ll 1$.
In the simple case of a linear drift function, i.e., when $a_i(\mathbf{x}(t)) = \mu _i x_i$ with $\mu_i < 0$, the probability density function of $x_i$ approaches a Gaussian with
(finite) variance $\mu_i$.
The time constant of the approach of the variance to equilibrium is $-1/\mu_i$ for $i=1,\ldots,m$ and $-\epsilon/\mu_i$ for $i=m+1,\ldots,n$ \cite{lelievre2013optimal}.
Thus, the last $n-m$ variables of \eqref{eq:general_SDE} rapidly approach a local equilibrium measure and exhibit fast dynamics, while the first $m$ variables exhibit slow dynamics.
A short burst of simulation will yield a cloud of points which is broadly distributed in the fast directions but narrowly distributed in the slow ones.
With the appropriate conditions on $a_i(\mathbf{x})$, the same can be said for more general drift
functions, where $\mu_i$ are the eigenvalues of the Jacobian of $\mathbf{a}(\mathbf{x}) = \begin{bmatrix} a_1(\mathbf{x}) & \cdots & a_n(\mathbf{x}) \end{bmatrix}^T$ \cite{villani2009hypocoercivity}.
Therefore, \eqref{eq:general_SDE} defines an $n$-dimensional stochastic system with $m$ slow
variables and $n-m$ fast variables, and $\epsilon$ defines the time scale separation.
The ratio of the powers of $\epsilon$ in the drift and diffusion terms in \eqref{eq:general_SDE} is essential,
as we require the square of the diffusivity to be of the same order as the drift as $\epsilon \rightarrow 0$ \cite{berglund2003geometric}.
If the diffusivity is larger, then, as $\epsilon \rightarrow 0$, the equilibrium measure will be
unbounded.
Conversely, if the diffusivity is smaller, the equilibrium measure will go to $0$ as $\epsilon \rightarrow 0$.

Assuming the sample average of $a_i(\mathbf{x})$ converges to a distribution which is only a function of the slow variables, then by the averaging principle \cite{freidlin2012random}, we can write a reduced SDE in {\em only} the slow variables $x_1, \dots, x_m$.
The aim of our work is to show how we can detect such slow variables {\em automatically} from data, in order to help inform modeling efforts and aid in the writing of such reduced stochastic models.
In general, we are not given the variables $\mathbf{x}(t)$ from the original SDE system, but instead, we are given some {\em observations}
 in the form $\mathbf{y}(t) = \mathbf{f} (\mathbf{x}(t))$.
We assume that $\mathbf{f}: \mathbb{R}^n \mapsto \mathbb{R}^d$, $n \le d$, is a deterministic (possibly nonlinear) function whose image is an $n$-dimensional manifold $\mathcal{M}$ in $\mathbb{R}^d$.
%
%
For our analysis, we require $\mathbf{g} = \mathbf{f} ^{-1}$ to be well-defined on $\mathcal{M}$, and both $\mathbf{f}$ and $\mathbf{g}$ to be continuously differentiable to fourth order.
Given data $\mathbf{y}(t_1),\ldots,\mathbf{y}(t_N)$ on $\mathcal{M}$ we would like to recover a parametrization of the data that is one-to-one with the
slow variables $x_1, \dots, x_m$.

\section{Local Invariant Metrics}

In order to recover the slow variables from data, we will utilize a local metric that collapses the fast directions.
Typically, such a metric averages out the fast variables.
However, simple averages are inadequate to describe data which is observed through a complicated nonlinear function.
Instead, we propose to use the Mahalanobis distance, which measures distances normalized by the respective variances in each local principal direction.
Using this metric, we still retain information about both the fast and slow directions and can
more clearly observe complex dynamic behavior within the data set.

If two points $\mathbf{x}(t_1)$ and $\mathbf{x}(t_2)$ are drawn from an $n$-dimensional
Gaussian distribution with covariance $\mathbf{C}_x$, the Mahalanobis distance between the points is defined as \cite{mahalanobis1936generalized}
\begin{equation}
	\| \mathbf{x}(t_1) - \mathbf{x}(t_2) \| _M = \sqrt{ (\mathbf{x}(t_1) - \mathbf{x}(t_2))^T \mathbf{C}_x^{-1} (\mathbf{x}(t_1) - \mathbf{x}(t_2) )  }.
\end{equation}
In particular, 
for \eqref{eq:general_SDE}, whose covariance does not depend on $\mathbf{x}$,  $\mathbf{C}_x^{-1} = \mathrm{diag}(e_1, \ldots, e_n)$ is a constant matrix where
\begin{equation} \label{eq:e_def}
\begin{aligned}
e_i =& 1, \: & 1 \le i \le m \\
e_i =& \epsilon, \: & m+1 \le i \le n,
\end{aligned}
\end{equation}
and the Mahalanobis distance between samples is
\begin{equation} \label{eq:rescale_x_dist}
\| \mathbf{x}(t_2) - \mathbf{x}(t_1) \|^2_M = \sum_{i=1}^n e_i \left( x_i(t_2) - x_i(t_1) \right)^2.
\end{equation}
Note that in \eqref{eq:rescale_x_dist}, the fast variables are collapsed and become $\mathcal{O}(\sqrt{\epsilon})$ small,
and so this metric is implicitly insensitive to variations in the fast variables.
The metric \eqref{eq:rescale_x_dist} can be rewritten as
\begin{equation} \label{eq:norm_z}
\| \mathbf{x}(t_2) - \mathbf{x}(t_1) \|^2_M = \| \mathbf{z}(t_2) - \mathbf{z}(t_1) \|^2_2
\end{equation}
where
\begin{equation} \label{eq:general_rescale}
z_i(t) = \sqrt{e_i} x_i(t).
\end{equation}
$\mathbf{z}(t)$ is a stochastic process of the same dimension as $\mathbf{x}(t)$, rescaled so that each variable has unit diffusivity.
This rescaling transforms our problem from one of detecting the slow variables within dynamic data to one of traditional data mining.
The Mahalanobis distance incorporates information about the dynamics and relevant time scales, so that using traditional data mining techniques with this metric will allow us to detect the slow variables in our data \cite{singer2009detecting}.
It is important to note that, in practice, we {\em never construct} $\mathbf{z}(t)$ {\em explicitly}.
It was shown in \cite{singer2008non} that, assuming $\mathbf{f}$ is bilipschitz, the Mahalanobis distance can be extended to approximate (to fourth order) the Euclidean distance between the rescaled samples $\mathbf{z}(t)$ from accessible $\mathbf{y}(t) = \mathbf{f} (\mathbf{x}(t))$,
\begin{equation} \label{eq:mahalanobis}
\| \mathbf{y}(t_2) - \mathbf{y}(t_1) \|^2_M = \| \mathbf{z}(t_2) - \mathbf{z}(t_1) \|^2_2 + \mathcal{O}(\| \mathbf{y}(t_2) - \mathbf{y}(t_1) \|^4_2).
\end{equation}
This approximation is accurate when $\| \mathbf{y}(t_2) - \mathbf{y}(t_1) \|$ is small.
Because we will integrate these distances into a manifold learning algorithm which only considers local distances, we can recover a parametrization of the data which is consistent with the underlying system variables $\mathbf{x}(t)$, even when the data are obscured by a function $\mathbf{f}$.
In Section~\ref{sec:analysis}, we will show how we can approximate this distance directly from data $\mathbf{y}(t)$.

\section{Diffusion Maps for Global Parametrization}

From pairwise distances, we want to extract a {\em global} parametrization of the data that represents the slow variables.
We will use diffusion maps \cite{Coifman2006, coifman2005geometric}, a kernel-based manifold learning technique, to extract a global parametrization using the local distances that we described in the previous section.
Given data $\mathbf{y}(t_1), \dots, \mathbf{y}(t_N)$, we first construct the kernel matrix $\mathbf{W} \in \mathbb{R}^{N \times N}$, where
\begin{equation} \label{eq:dmaps_kernel}
W_{ij} = \exp \left( -\frac{\|\mathbf{y}(t_i) - \mathbf{y}(t_j) \|^2}{\sigma_{kernel}^2} \right).
\end{equation}
Here, $\| \cdot \|$ denotes the appropriate norm (in our case, the Mahalanobis distance), and $\sigma_{kernel}$ is the kernel scale
and denotes a characteristic distance within the data set.
Note that $\sigma_{kernel}$ induces a notion of locality: if $\|\mathbf{y}(t_i) - \mathbf{y}(t_j) \| \gg \sigma_{kernel}$, then $W_{ij}$ is negligible.
Therefore, we only need our metric to be informative within a ball of radius $\sigma _{kernel}$.
We then construct the diagonal matrix $\mathbf{D} \in \mathbb{R}^{N \times N}$, with
\begin{equation}
D_{ii} = \sum_{j=1}^N W_{ij}.
\end{equation}
We compute the eigenvalues $\lambda_0, \dots, \lambda_{N-1}$ and eigenvectors $\phi_0, \dots, \phi_{N-1}$ of the matrix $\mathbf{A} = \mathbf{D}^{-1}\mathbf{W}$, and order them such that $1 = \lambda_0 \ge |\lambda_1| \ge \dots \ge |\lambda_{N-1}|$.
$\phi_0 = \begin{bmatrix} 1 & 1 & \cdots & 1 \end{bmatrix}^T$
is the trivial eigenvector; the next few eigenvectors provide embedding coordinates for the data, so that $\phi_j(i)$, the $i^{th}$ entry of $\phi_j$, provides the $j^{th}$ embedding coordinate for $\mathbf{y}(t_i)$ (modulo higher harmonics which characterize
the same direction in the data \cite{ferguson2010systematic}).

\section{Estimation of the Mahalanobis Distance} \label{sec:analysis}

\begin{figure}[t]
\centering
\epsfig{width=\textwidth, file=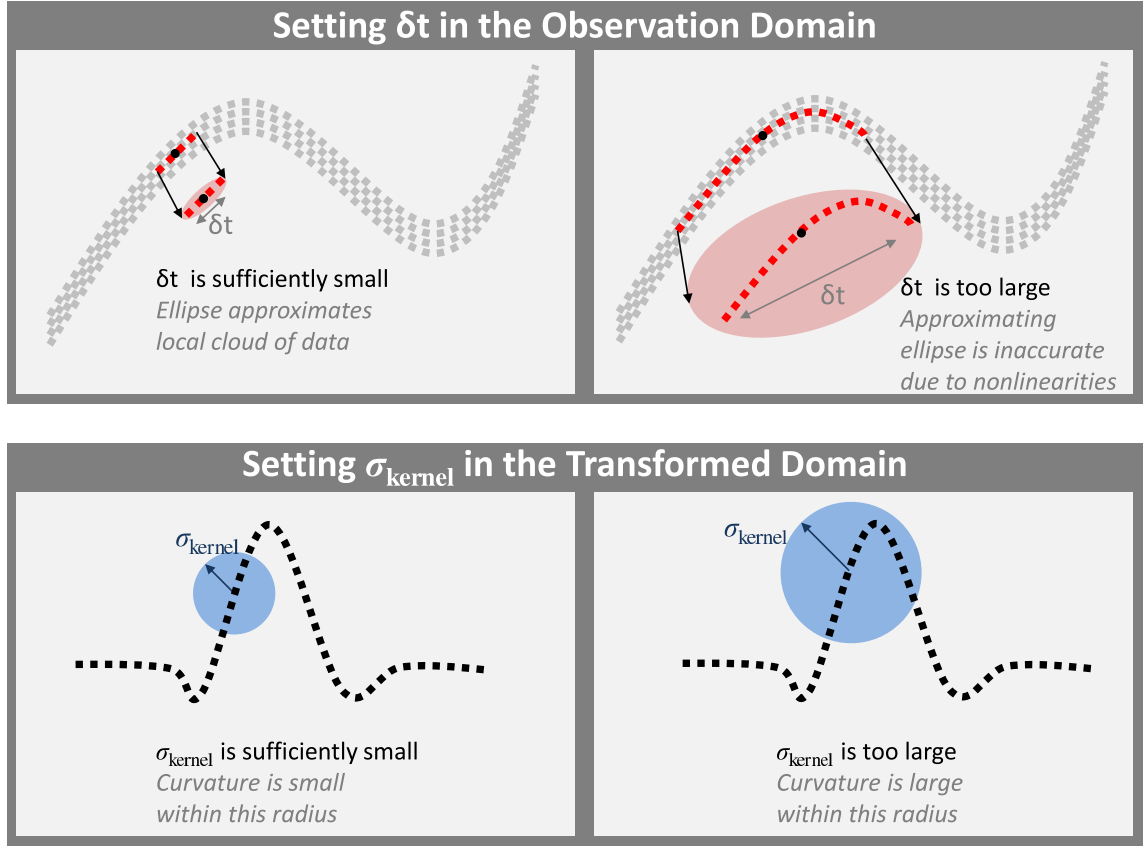}
\caption{Illustration of how to choose $\delta t$ and $\sigma_{kernel}$ appropriately. 	Curvature effects and other nonlinearities should be negligible within a time window $\delta t$ and within a ball of radius $\sigma_{kernel}$.}
\label{fig:schematic}
\end{figure}

As previously mentioned, we do not have access to the original variables $\mathbf{x}(t)$ from the underlying original SDE system.
Instead, we only have measurements $\mathbf{y}(t) = \mathbf{f} (\mathbf{x}(t))$, and we want to estimate the Mahalanobis distance between the $\mathbf{x}$ variables from observations $\mathbf{y}(t)$.
The traditional Mahalanobis distance is defined for a fixed distribution,
whereas here we are dealing with a distribution that possibly changes as a function of position due to
nonlinearities in the observation function $\mathbf{f}$ and in the drift $\mathbf{a}(\mathbf{x})$.
Consequently, we use the following modified definition for the Mahalanobis distance between two points,
\begin{equation} \label{eq:mahalanobis_distance}
 \| \mathbf{y}(t_2) - \mathbf{y}(t_1) \|^2_M =
 \frac{1}{2} (\mathbf{y}(t_2) - \mathbf{y}(t_1))^T \left( \mathbf{C}^{\dagger}(\mathbf{y}(t_1)) + \mathbf{C}^{\dagger}(\mathbf{y}(t_2)) \right) (\mathbf{y}(t_2) - \mathbf{y}(t_1)),
 \end{equation}
where $\mathbf{C}(\mathbf{y}(t))$ is the covariance of the observed stochastic process {\em at the point} $\mathbf{y}(t)$,
and $\dagger$ denotes the Moore-Penrose pseudoinverse (since $d$ may exceed $n$).

To motivate this definition of the Mahalanobis distance, we first consider the simple linear case where $\mathbf{f} (\mathbf{x}) = \mathbf{A} \mathbf{x}$, with $\mathbf{A} \in \mathbb{R}^{d \times n}$.
The covariance of the observed stochastic process $\mathbf{f} (\mathbf{x})$ is given by $\mathbf{C}=\mathbf{AC}_x\mathbf{A}^T$.
Let $\mathbf{A}= \mathbf{U} \mathbf{\Lambda} \mathbf{V}^T$ be the singular value decomposition (SVD) of $\mathbf{A}$, where $\mathbf{U} \in \mathbb{R}^{d \times n}$, $\mathbf{\Lambda} \in \mathbb{R}^{n \times n}$, and $\mathbf{V} \in \mathbb{R}^{n \times n}$.
The pseudoinverse of the covariance matrix is $\mathbf{C}^{\dagger} = \mathbf{U} \mathbf{\Lambda}^{-1} \mathbf{V}^T \mathbf{C}_x^{-1} \mathbf{V \Lambda} ^{-1} \mathbf{U}^T$. Consequently, the Mahalanobis distance \eqref{eq:mahalanobis_distance} is reduced to
\begin{equation}
\begin{aligned}
\| \mathbf{y}(t_2) - \mathbf{y}(t_1) \|^2_M &= (\mathbf{y}(t_2) - \mathbf{y}(t_1))^T  \mathbf{C}^{\dagger} (\mathbf{y}(t_2) - \mathbf{y}(t_1)) \\
 &= (\mathbf{x}(t_2) - \mathbf{x}(t_1))^T \mathbf{A}^T \mathbf{C}_x^{-1} \mathbf{A} (\mathbf{x}(t_2) - \mathbf{x}(t_1)) \\
 &= (\mathbf{x}(t_2) - \mathbf{x}(t_1))^T \mathbf{V \Lambda U}^T \mathbf{ U \Lambda}^{-1} \mathbf{V}^T \mathbf{C}_x^{-1} \mathbf{V \Lambda} ^{-1} \mathbf{U}^T \mathbf{U \Lambda V}^T (\mathbf{x}(t_2) - \mathbf{x}(t_1)) \\
 &= (\mathbf{x}(t_2) - \mathbf{x}(t_1))^T \mathbf{C}_x^{-1}  (\mathbf{x}(t_2) - \mathbf{x}(t_1)) \\
 &= \| \mathbf{x}(t_2) - \mathbf{x}(t_1) \|^2_M =  \| \mathbf{z}(t_2) - \mathbf{z}(t_1) \|^2_2 .
\end{aligned}
\end{equation}
Hence evaluating the Mahalanobis distances of the observations $\mathbf{y}(t) = \mathbf{f}(\mathbf{x}(t))$ using \eqref{eq:mahalanobis_distance} allows us to estimate the Euclidean distances of the rescaled variables $\mathbf{z}$ (in which the fast coordinates are collapsed).

Following \cite{singer2008non}, we will show via Taylor expansion that the Mahalanobis distance between the observations \eqref{eq:mahalanobis_distance} approximates the Euclidean
distance {\em in the rescaled variables} for general nonlinear observation functions $\mathbf{f}$ (provided $\mathbf{f}$ is bilipschitz and both $\mathbf{f}$ and $\mathbf{f}^{-1}$ are differentiable to fourth order).
%
\eqref{eq:mahalanobis_distance} cannot be evaluated directly since we do not have access to the covariance matrices, so we will instead estimate the covariances directly from data.
We can estimate the covariance $\mathbf{C}(\mathbf{y}(t_0))$ empirically from a set of values $\mathbf{y}(t_1), \dots, \mathbf{y}(t_q)$ drawn from the local distribution at $\mathbf{y}(t_0)$.
One way to obtain such a set of points is to run $q$ simulations for a short time, $\delta t$, each starting from $\mathbf{y}(t_0)$.
Alternatively, we can consider a single time series of length $q \delta t$ starting from $\mathbf{y}(t_0)$, and then estimate the covariance
from the increments $\Delta \mathbf{y}(t_i) = \mathbf{y}(t_i) -\mathbf{y}(t_{i-1})$.
Although we will present analysis and results for the first type of estimation, the second case is often more practical in practice.
%

Errors in our estimation of the Mahalanobis distance arise from three sources.
One source of error is approximating the function $\mathbf{f}$ locally as a linear function by truncating the Taylor expansion of $\mathbf{f}$ at first order.
An additional source of error arises from disregarding the drift in the stochastic process, and assuming that samples are drawn from a Gaussian distribution.
The third source comes from finite sampling effects.
In this work, we will address and discuss the first two sources of error (the finite sampling effects are the subject of future research).
We can control the effects of the errors due to truncation of the Taylor expansion by adjusting $\sigma_{kernel}$; the higher-order terms in this expansion will be small for points which are close, such that adjusting $\sigma_{kernel}$ will allow us to only consider distances which are
sufficiently accurate in our overall computation scheme.
Furthermore, we can control the errors incurred by disregarding the drift by adjusting the time scale of our simulation bursts $\delta t$.
Figure~\ref{fig:schematic} illustrates some of the issues in choosing the sizes $\delta t$ (or $q \delta t$ if the alternate method is used) and the parameter $\sigma_{kernel}$.
We will present both analytical results for the error bounds, as well as an empirical methodology to set the parameters $\sigma_{kernel}$ and $\delta t$ for our method to accurately recover the slow variable(s).

\subsection{Error due to the observation function $\mathbf{f}$}

We want to relate the distance in the rescaled space, $\|\mathbf{z}(t_2) - \mathbf{z}(t_1)\|_2$, to the estimated Mahalanobis distance between the observations $\| \mathbf{y}(t_2) - \mathbf{y}(t_1)\|_M$.
We define the error incurred by using the Mahalanobis distance to approximate the true distance as
\begin{equation}
E_M(\mathbf{y}(t_1), \mathbf{y}(t_2)) = \|\mathbf{z}(t_2) - \mathbf{z}(t_1)\|_2^2 - \| \mathbf{y}(t_2) - \mathbf{y}(t_1)\|^2_M .
\end{equation}
By Taylor expansion of $\mathbf{g}(y) = \mathbf{f}^{-1}(y)$ around $\mathbf{y}(t_1)$ and $\mathbf{y}(t_2)$ and averaging the two expansions, we obtain
\begin{equation} \label{eq:mahanaobis_error}
\begin{aligned}
E_M\left( \mathbf{y}(t_1), \mathbf{y}(t_2) \right)
 =&
\begin{aligned}[t]
 \frac{1}{2} \sum_{i=1}^n \sum_{jkl=1}^{d} &
\left( g_{i, (j)} (\mathbf{y}(t_1)) g_{i, (k,l)} (\mathbf{y}(t_1)) -  g_{i, (j)} (\mathbf{y}(t_2)) g_{i, (k,l)} (\mathbf{y}(t_2)) \right) \times \\
& ({y}_j(t_2) - {y}_j(t_1)) ({y}_k(t_2) - {y}_k(t_1))({y}_l(t_2) - {y}_l(t_1))
\end{aligned} \\
+&
\begin{aligned}[t]
\frac{1}{8} \sum_{i=1}^n \sum_{jklm=1}^d  &
\left( g_{i, (j,k)} (\mathbf{y}(t_1)) g_{i, (l,m)} (\mathbf{y}(t_1)) +  g_{i, (j,k)} (\mathbf{y}(t_2)) g_{i, (l,m)} (\mathbf{y}(t_2)) \right) \times
 \\
&({y}_j(t_2) - {y}_j(t_1))  ({y}_k(t_2) - {y}_k(t_1))({y}_l(t_2) - {y}_l(t_1)) ({y}_m(t_2) - {y}_m(t_1))
\end{aligned} \\
+&
\begin{aligned} [t]
\frac{1}{6} \sum_{i=1}^n \sum_{jklm=1}^d &
\left( g_{i, (j)} (\mathbf{y}(t_1)) g_{i, (k,l,m)} (\mathbf{y}(t_1)) +  g_{i, (j)} (\mathbf{y}(t_2)) g_{i, (k,l,m)} (\mathbf{y}(t_2)) \right) \times \\
& ({y}_j(t_2) - {y}_j(t_1))  ({y}_k(t_2) - {y}_k(t_1))({y}_l(t_2) - {y}_l(t_1))({y}_m(t_2) - {y}_m(t_1))
\end{aligned} \\
+& \mathcal{O} \left(\| \mathbf{y}(t_2) - \mathbf{y}(t_1) \|^6_2 \right) ,
\end{aligned}
\end{equation}
where
\begin{equation}
\begin{aligned}
g_{i,(j)} &= \sqrt{e_i} \frac{\partial g_i}{\partial y_j}
\\
g_{i,(j,k)} &= \sqrt{e_i}  \frac{\partial^2 g_i}{\partial y_j \partial y_k}
\\
g_{i,(j,k,l)} &= \sqrt{e_i}  \frac{\partial^3 g_i}{\partial y_j \partial y_k \partial y_l} .
\end{aligned}
\end{equation}
In \cite{singer2008non}, it was shown that the error incurred by using the Mahalanobis distance to approximate the $L_2$-distance between points $\mathbf{z}(t)$ is $\mathcal{O} (\|\mathbf{y}_1 - \mathbf{y}_2 \|_2^4 )$ (see the Supplementary Materials for details).
We now see from \eqref{eq:mahanaobis_error} that the error is an explicit function of the second- and higher-order derivatives of $\mathbf{g} = \mathbf{f}^{-1}$ and the distance between samples $\| \mathbf{y}(t_2) - \mathbf{y}(t_1) \|_2$.
We would like to note that this error does not depend on the dynamics of the underlying stochastic process (as we assume the covariances at each point on the manifold are known), but is only a function of the measurement function $\mathbf{f}$.
The parameter $\sigma_{kernel}$ in the diffusion maps calculation determines how much $E_M$ contributes to the overall analysis.
From \eqref{eq:dmaps_kernel}, distances which are much greater than $\sigma_{kernel}$ are negligible in the diffusion maps computation because of the exponential kernel.
Therefore, we want to choose $\sigma_{kernel}^2$ on the order of $\|\mathbf{y}(t_2) - \mathbf{y}(t_1)\|^2_M$ in a regime where $| E_M(\mathbf{y}(t_1), \mathbf{y}(t_2))|  \ll \|\mathbf{y}(t_2) - \mathbf{y}(t_1)\|^2_M$.
This is illustrated in Figure~\ref{fig:schematic}, where we want to choose $\sigma_{kernel}$ small enough so that the curvature and other nonlinear effects (captured in the error term $E_M$) are negligible.
This will ensure that the errors in the Mahalanobis distance approximation do not greatly effect our overall analysis.

On first inspection, it would appear that our analysis indicates that $\sigma_{kernel}$ should be chosen arbitrarily small.
However, to obtain a meaningful parametrization of the data set, there must be a nonnegligible number of data points within a ball of radius $\sigma_{kernel}$ around each sample.
Therefore, the sampling density on the underlying manifold provides a lower bound for $\sigma_{kernel}$.

\subsection{Error due to the dynamics} \label{subsec:cov_est}

To compute the Mahalanobis distance in \eqref{eq:mahalanobis_distance}, we require $\mathbf{C}$, the covariance of the observed stochastic process $\mathbf{y}(t) = \mathbf{f}( \mathbf{x}(t))$.
We will use simulation bursts to locally explore the dynamics on the manifold of observations in order to estimate the covariance at a point $\mathbf{y}(t)$ from data \cite{talmon2014manifold, talmon2014intrinsic}.
We write the elements of the estimated covariance $\hat{\mathbf{C}}(\mathbf{y}(t), \delta t)$ as
\begin{equation}\label{eq:estimated_cov_expected_value}
\hat{C}_{ij}(\mathbf{y}(t), \delta t)
=
\frac{1}{\delta t} \left( \mathbb{E} \left[ y_i (t+\delta t) y_j (t+ \delta t) \mid \mathbf{y}(t) \right]
- \mathbb{E} \left[ y_i (t+\delta t) \mid \mathbf{y}(t) \right] \mathbb{E} \left[ y_j (t+\delta t) \mid \mathbf{y}(t) \right] \right) ,
\end{equation}
where $\delta t > 0$ is the length of the simulation burst.

Due to the drift in the stochastic process and the (perhaps nonlinear) measurement function $\mathbf{f}$, we incur some error by approximating the covariance at a point $\mathbf{y}(t)$ using simulations of length $\delta t > 0$.
Define the error in this approximation as
\begin{equation}
\mathbf{E}_C(\mathbf{y}(t), \delta t) = \hat{\mathbf{C}}(\mathbf{y}(t), \delta t) - \mathbf{C}(\mathbf{y}(t)).
\end{equation}
%
%
By It\^{o}-Taylor expansion of $\mathbf{f}$ and $\mathbf{x}(t)$ \cite{kloeden1992numerical},
\begin{equation} \label{eq:cov_error}
\begin{aligned}
E_{C, ij} (\mathbf{x}(t), \delta t) =&
 \frac{1}{\delta t} \sum_{k=1}^n f_{i,(k)}(\mathbf{x}(t)) \mathbb{E} \left[ \int_t^{t+\delta t} \left( \int_{s_2}^{t+\delta t} f_{j,(k,0)}(\mathbf{x}(s_1)) ds_1
+ \int_t^{s_2} f_{j,(0,k)}(\mathbf{x}(s_1)) ds_1 \right) ds_2 \right] \\
&+  \frac{1}{\delta t} \sum_{k=1}^n f_{j,(k)}(\mathbf{x}(t))  \mathbb{E} \left[ \int_t^{t+\delta t} \left( \int_{s_2}^{t + \delta t} f_{i,(k,0)}(\mathbf{x}(s_1)) ds_1
+  \int_t^{s_2} f_{i,(0,k)}(\mathbf{x}(s_1)) ds_1 \right) ds_2 \right] \\
&+  \frac{1}{\delta t} \sum_{k,l=1}^n \mathbb{E} \left[ \int_t^{t+\delta t}\left( \int_t^{s_2} f_{i,(k,l)}(\mathbf{x}(s_1)) dW_{s_1, k}  \right) \left(  \int_t^{s_2} f_{j,(k,l)}(\mathbf{x}(s_1)) dW_{s_1, k} \right) ds_2 \right] \\
&+ \mathcal{O} (\delta t^{3/2})
\end{aligned}
\end{equation}
where
\begin{equation}
\begin{aligned}
f_{i,(k)} &= \frac{1}{\sqrt{e_k}} \frac{\partial f_i}{\partial x_k}
\\
f_{i,(k,l)} &= \frac{1}{\sqrt{e_k e_l}} \frac{\partial^2 f_i}{\partial x_k \partial x_l}
\\
f_{i,(k,0)} &= \frac{1}{\sqrt{e_k}} \sum_{l=1}^n \left( \frac{\partial}{\partial x_k} \left( \frac{a_l(\mathbf{x})}{e_l} \frac{\partial f_i}{\partial x_l} \right) + \frac{1}{2 e_l} \frac{\partial^3 f_i}{\partial x_k \partial x_l^2} \right)
\\
f_{i,(0, k)} &= \frac{1}{\sqrt{e_k}} \sum_{l=1}^n \left( \frac{a_l(\mathbf{x})}{e_l} \frac{\partial^2 f_i}{\partial x_k \partial x_l} +\frac{1}{2 e_l}  \frac{\partial^3 f_i}{\partial x_k \partial^2 x_l} \right).
\end{aligned}
\end{equation}
From \eqref{eq:cov_error}, the error in the covariance is $\mathcal{O}(\delta t)$ (as the $ds$ integrals are each $\mathcal{O}(\delta t)$ and the $dW$ integrals are each $\mathcal{O}(\sqrt{\delta t})$) and a function of the derivatives of the observation function $\mathbf{f}$ and the drift $\mathbf{a}$.
We want to set $\delta t$ such that $\|\mathbf{E}_C \| \ll \| \mathbf{C} \|$
(this is illustrated in Figure~\ref{fig:schematic}), so that the estimated covariances are accurate.
Note that in practice, we compute $\hat{\mathbf{C}}$ by running many simulations of length $\delta t$ starting from $\mathbf{x}(t)$, and use the sample average to approximate the expected values in \eqref{eq:estimated_cov_expected_value}.
We therefore incur additional error due to finite sampling; this error is ignored for the purposes of this analysis, and quantifying this error is the subject of future research.

Our analysis reveals that the errors decrease with decreasing $\delta t$; at first inspection, one would want to set $\delta t$ arbitrarily small to obtain the highest accuracy possible.
However, often in practice, one cannot obtain an arbitrarily refined sampling rate, such that a smaller $\delta t$ results in fewer samples with which to approximate the local covariance.%
When also accounting for these finite sampling errors, and one should take $\delta t$ as long as possible while still maintaining negligable errors from the observation function $\mathbf{f}$ and the drift $\mathbf{a}$.

\section{Illustrative Examples}

\begin{figure}[t]
\centering
\epsfig{width=0.5\textwidth, file=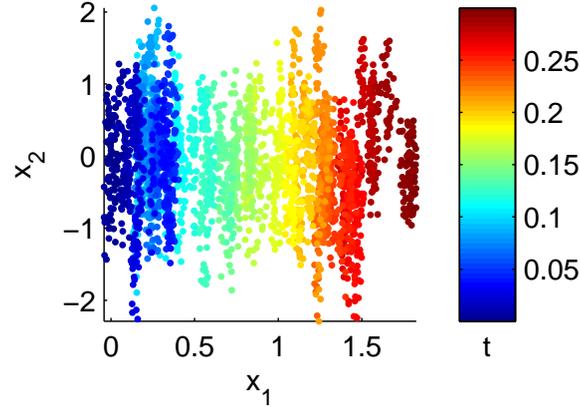}
\caption{Data, simulated from \eqref{eq:specific_SDE} with $a=3$ and $\epsilon = 10^{-3}$, for $3000$ time steps with $dt = 10^{-4}$. The data are colored by time.}
\label{fig:initial_data}
\end{figure}

For illustrative purposes, we consider the following two-dimensional SDE
\begin{equation} \label{eq:specific_SDE}
\begin{aligned}
dx_1(t) &=& adt &+& dW_1(t)\\
dx_2(t) &=& -\frac{x_2(t)}{\epsilon} dt &+& \frac{1}{\sqrt{\epsilon}} dW_2(t)
\end{aligned}
\end{equation}
where $a$ is an $\mathcal{O}(1)$ constant, as a specific example of \eqref{eq:general_SDE}.
$x_1$ is the slow variable, and $x_2$ is a fast noise whose equilibrium measure is bounded and $\mathcal{O}(1)$.
Figure~\ref{fig:initial_data} shows data simulated from this SDE colored by time.
We would like to recover a parametrization of this data which is one-to-one with the slow variable $x_1$.

\subsection{Linear function} \label{subsec:linear_example}

In the first example, our observation function $\mathbf{f}$ will be the identity function,
\begin{equation} \label{eq:linear_transform}
\begin{aligned}
\begin{bmatrix}
y_1(t) \\ y_2(t)
\end{bmatrix} &=&
\mathbf{f}(\mathbf{x}(t)) &=&
\begin{bmatrix} x_1(t) \\ x_2(t) \end{bmatrix} \\
\mathbf{g}(\mathbf{y}(t)) &=& \mathbf{f}^{-1} (\mathbf{y}(t)) &=& \begin{bmatrix} y_1(t) \\ y_2(t) \end{bmatrix}
\end{aligned}
\end{equation}
where the fast and slow variables remain uncoupled.
In this case, there is no error incurred due to the measurement function $\mathbf{f}$ ($E_M = 0$), as the second- and higher-order derivatives of $\mathbf{g}$ are identically 0.

\subsubsection{Importance of using the Mahalanobis distance}

\begin{figure}[t]
\centering
\begin{subfigure}{0.4\textwidth}
\centering
\epsfig{height=2in, file=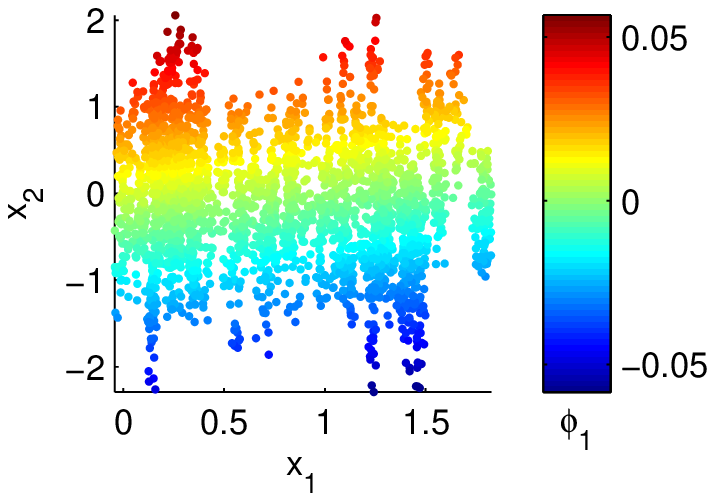}
\caption{}
\end{subfigure}
\begin{subfigure}{0.4\textwidth}
\centering
\epsfig{height=2in, file=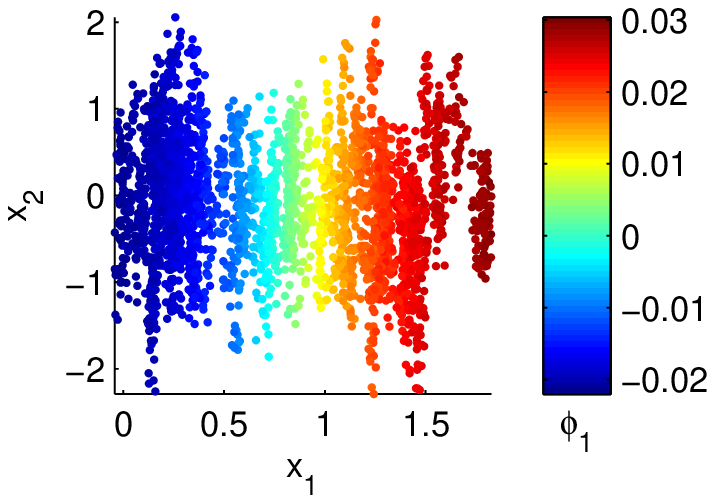}
\caption{}
\end{subfigure}
\caption{Comparison of using the Euclidean distance and the Mahalanobis distance in multiscale data mining. (a) The data from Figure~\ref{fig:initial_data}, colored by the first diffusion map coordinate when using the Euclidean distance in the kernel in \eqref{eq:dmaps_kernel}. Note that we do {\em not} recover the slow variable. (b) The data from Figure~\ref{fig:initial_data}, colored by the first diffusion map coordinate when using the Mahalanobis distance in the kernel in \eqref{eq:dmaps_kernel}. The good correspondence between this coordinate and the slow variable is visually obvious.}
\label{fig:NIV_versus_DMAPS}
\end{figure}

We want to demonstrate the utility of using the Mahalanobis distance compared to the typical Euclidean distance.
We compute the diffusion map embedding for the data in Figure~\ref{fig:initial_data},
using both the standard Euclidean distance and the Mahalanobis distance for the computation of the kernel in \eqref{eq:dmaps_kernel}.
The data, colored by $\phi_1$ using the two different metrics, are shown in Figure~\ref{fig:NIV_versus_DMAPS}.
When using the standard Euclidean distance which does not account for the underlying dynamics, the first diffusion maps recovers the fast variable $x_2$, suggesting the fast modes is the dominant scale purely in terms of data analysis (Figure~\ref{fig:NIV_versus_DMAPS}(a)).
In contrast, the slow variable is recovered when using the Mahalanobis distance, as the coloring in Figure~\ref{fig:NIV_versus_DMAPS}(b) (where the data are colored by the first diffusion maps variable) is consistent with the coloring in Figure~\ref{fig:initial_data} (where the data are colored by time).

\subsubsection{Errors in covariance estimation}

\begin{figure}[t]
\centering
\begin{subfigure}{0.4\textwidth}
\centering
\epsfig{height=2in, file=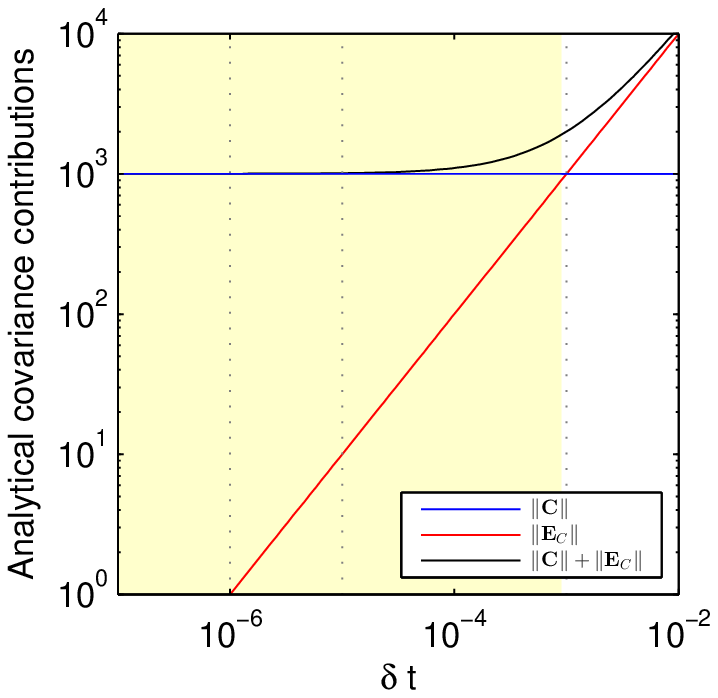}
\caption{}
\end{subfigure}
\begin{subfigure}{0.4\textwidth}
\centering
\epsfig{height=2in, file=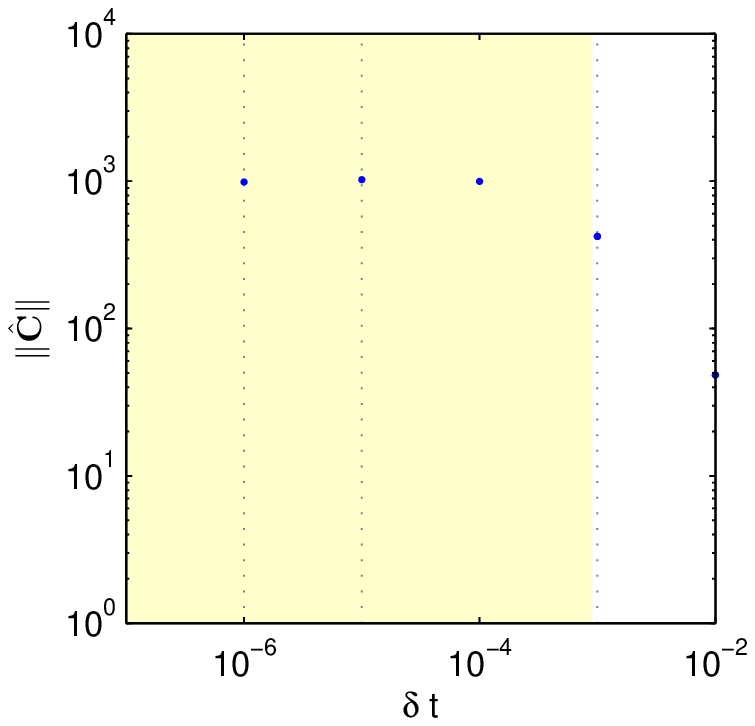}
\caption{}
\end{subfigure}
\caption{Errors in covariance estimation for linear example. (a) The analytical contributions to the covariance for the example in Section~\ref{subsec:linear_example} as a function of $\delta t$. (b) The average estimated covariance $\| \hat{ \mathbf{C}} \|$ as a function of $\delta t$. The average is computed over $10$ data points and using $50$ sample points to estimate each covariance. The shaded yellow region indicates the range of $\delta t$ over which the errors in the estimated covariance are less than the norm of the covariance. }
\label{fig:cov_error}
\end{figure}

For the example in \eqref{eq:linear_transform}, the analytical covariance is
 \begin{equation} \label{eq:cov_linear_example}
\mathbf{C}(\mathbf{x}(t)) =
\begin{bmatrix}
1 & 0 \\
0 & \frac{1}{\epsilon}
\end{bmatrix}.
\end{equation}
From \eqref{eq:cov_error}, we find
\begin{equation}
\mathbf{E}_C(\mathbf{x}(t), \delta t) =
\begin{bmatrix}
0 & 0 \\
0 & -\frac{\delta t}{\epsilon^2}
\end{bmatrix}
+ \mathcal{O} (\delta t^{3/2}) .
\end{equation}
Therefore, $\| \mathbf{C} \| = \mathcal{O} \left( \frac{1}{\epsilon} \right)$ and $\|\mathbf{E}_C \| = \mathcal{O}\left(\frac{\delta t}{\epsilon^2} \right)$ (provided $\frac{1}{\epsilon^2} \gg \sqrt{\delta t}$; this will be discussed further in Section~\ref{subsec:fastvar}).
These terms are shown in Figure~\ref{fig:cov_error}(a) as a function of $\delta t$.
We want to choose $\delta t$ in a regime where $\| \mathbf{E}_C \| \ll \| \mathbf{C} \|$ (the yellow shaded region in Figure~\ref{fig:cov_error} indicates where $\| \mathbf{E}_C \| < \| \mathbf{C} \|$), so that the errors in the estimated covariance are small with respect to the covariance.

When we do not analytically know the functions $\mathbf{f}$ or $\mathbf{g}$, we can find such a regime empirically by
estimating the covariance for several values of $\delta t$.
This provides an estimate of $\hat{\mathbf{C}} = \mathbf{C} + \mathbf{E}_C$ as a function of $\delta t$.
From Figure~\ref{fig:cov_error}(a), we expect a ``knee" in the plot of $\| \hat{\mathbf{C}} \|$ versus $\delta t$ when $\| \mathbf{E}_C \|$ becomes larger than $\| \mathbf{C}\|$.
Figure~\ref{fig:cov_error}(b) shows the empirical $\| \hat{\mathbf{C}} \|$ as a function of $\delta t$ for the data in Figure~\ref{fig:initial_data}, and the knee in this curve is consistent with the intersection in Figure~\ref{fig:cov_error}(a).

\subsubsection{Recovery of the fast variable} \label{subsec:fastvar}

Note that, for the example in \eqref{eq:linear_transform}, $\mathbf{E}_C$ is a constant diagonal matrix.
Therefore, taking $\delta t$ too large will not lead to nonlinear effects or mixing of the fast and slow variables.
Rather, changing $\delta t$ will only affect the perceived ratio of the fast and slow timescales.

To see this behavior in our diffusion maps results, we must first discuss the interpretation of the diffusion maps eigenspectrum.
The diffusion maps eigenvectors provide embedding coordinates for the data, and the corresponding eigenvalues provide a measure of the importance of each coordinate.
However, some eigenvectors can be harmonics of previous eigenvectors; for example, for a data set parameterized by a variable $x$, both $\cos x$ and $\cos 2x$ will appear as diffusion maps eigenvectors (see \cite{ferguson2010systematic} for a more detailed discussion).
These harmonics do not capture any new direction within the data set, but do appear as additional eigenvector/eigenvalue pairs.
Therefore, for the two-dimensional data considered here, the fast variable will not necessarily appear as the second (non-trivial) eigenvector.
As the time scale separation increases, the relative importance of the slow and fast directions will also increase.
This implies that the eigenvalue corresponding to the eigenvector which parameterizes the fast direction will decrease, and the number of harmonics of the slow mode which appear before the fast mode will increase.

\begin{figure}[t]
\epsfig{width=0.3\textwidth, file=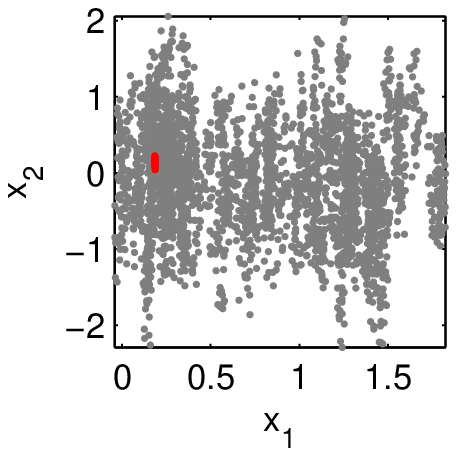}
\epsfig{width=0.3\textwidth, file=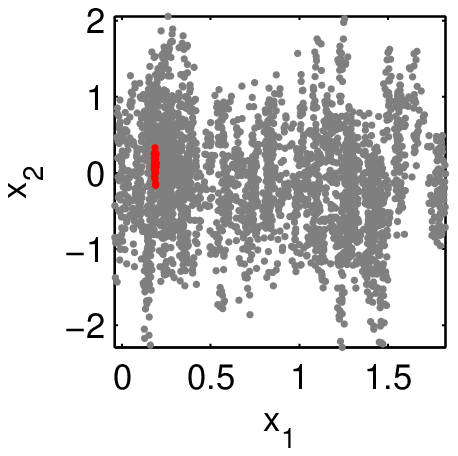}
\epsfig{width=0.3\textwidth, file=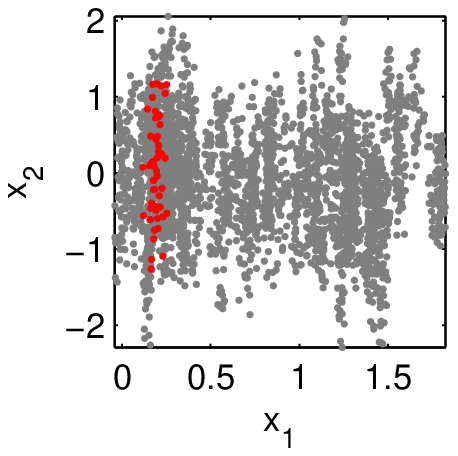}

\epsfig{width=0.3\textwidth, file=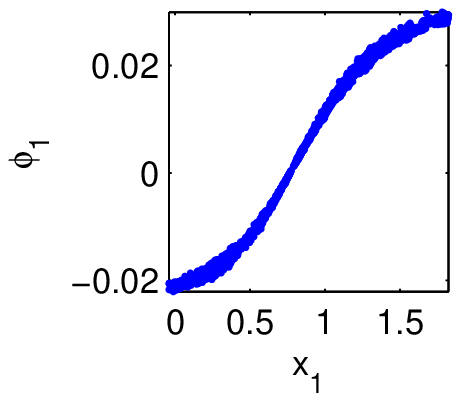}
\epsfig{width=0.3\textwidth, file=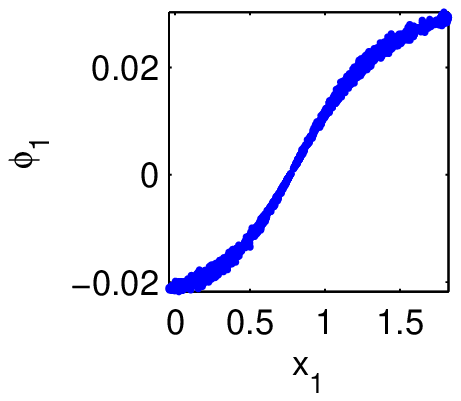}
\epsfig{width=0.3\textwidth, file=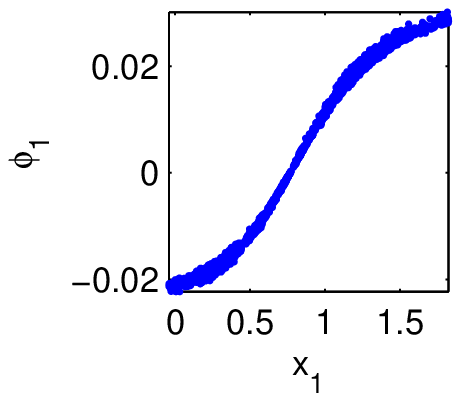}

\epsfig{width=0.3\textwidth, file=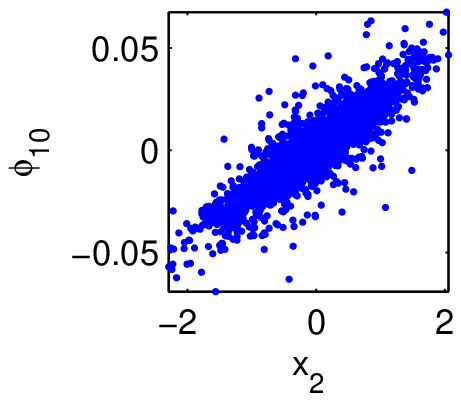}
\epsfig{width=0.3\textwidth, file=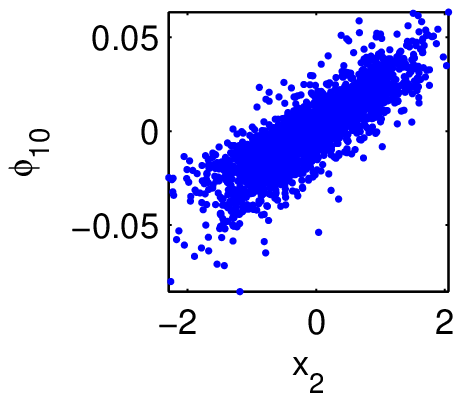}
\epsfig{width=0.3\textwidth, file=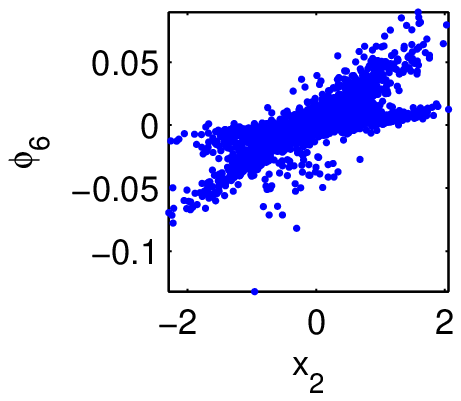}

\epsfig{width=0.3\textwidth, file=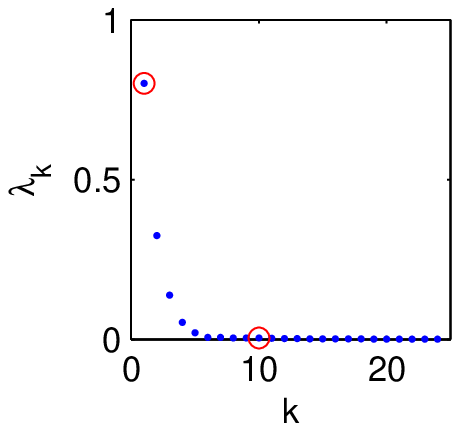}
\epsfig{width=0.3\textwidth, file=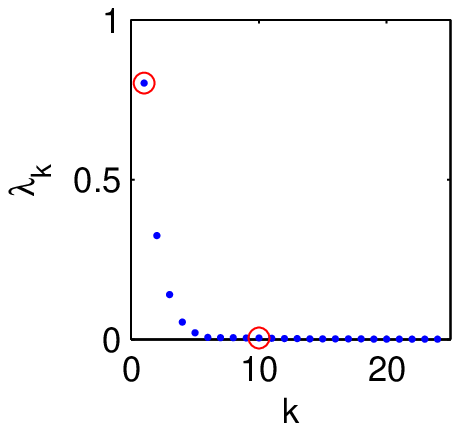}
\epsfig{width=0.3\textwidth, file=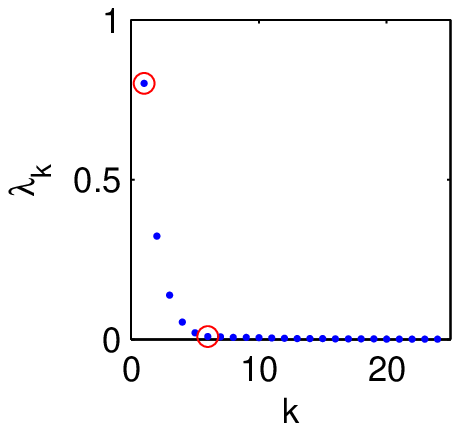}

\caption{Relationship between changing $\delta t$ and recovery of the variables. From left to right, the columns correspond to $\delta t = 10^{-6}, 10^{-5}, 10^{-3}$.  (Row~1) Data (gray) and representative burst (red) used to estimate the local covariance. (Row~2) Correlation between first diffusion maps coordinate and the slow variable $x_1$. (Row~3) Correlation between the relevant diffusion maps coordinate and the fast variable $x_2$. Note that for $\delta t = 10^{-6}$ and $\delta t = 10^{-5}$, $x_2$ is correlated with $\phi_{10}$. When $\delta t = 10^{-3}$, $x_2$ is correlated with $\phi_6$. (Row 4) Diffusion maps eigenvalue spectra. The eigenvalues corresponding to the coordinates for the slow and fast modes are indicated by red circles. Note that when $\delta t$ is too large, the apparent time scale separation decreases and the coordinate corresponding to the fast variable appears earlier in the spectrum. }
\label{fig:recover_fast}
\end{figure}

Figure~\ref{fig:recover_fast} shows results for three different values of $\delta t$ (the corresponding values are indicated by the dashed lines in Figure~\ref{fig:cov_error}).
When the time scale of the simulation burst used to estimate the local covariance (indicated by the red clouds in the top row of figures), is sufficiently shorter than that of the equilibration time of the fast variable, the estimated local covariance is accurate and the fast variable is collapsed significantly relative to the slow variable.
This means that the fast variable is recovered {\em very} far down in the diffusion maps eigenvectors.
The left two columns of Figure~\ref{fig:recover_fast} show that, for this example, when the simulation burst is shorter than the equilibration time, the fast variable is recovered as $\phi_{10}$.
However, if the time scale of the burst is {\em longer} than the saturation time of the fast variable, the estimated covariance changes: the variance in the slow direction continues to grow, while the variance in the fast direction is fixed.
This means that the apparent time scale separation is smaller, the collapse of the fast variable is less pronounced relative to the slow variable, and the fast variable is recovered in an earlier eigenvector (in our ordering of the spectrum).
The right column of Figure~\ref{fig:recover_fast} shows that, when the burst is now longer than the equilibration time, the fast variable appears earlier in the eigenvalue spectrum and is recovered as $\phi_6$.

\subsection{Nonlinear observation function} \label{subsec:nonlinear_example}

\begin{figure}[t]
\centering
\epsfig{width=0.5\textwidth, file=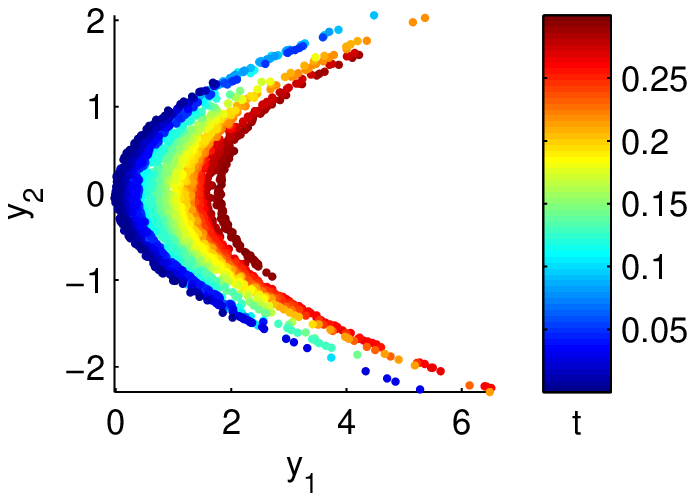}
\caption{The data from Figure~\ref{fig:initial_data}, transformed by $\mathbf{f}$ in \eqref{eq:nonlinear_function}}
\label{fig:initial_data_nonlinear}
\end{figure}

In the second example, our data will be warped into ``half-moon" shapes via the function
\begin{equation} \label{eq:nonlinear_function}
\begin{aligned}
\begin{bmatrix}
y_1(t) \\ y_2(t)
\end{bmatrix} &=&
\mathbf{f}(\mathbf{x}(t)) &=&
\begin{bmatrix}
x_1(t) + x_2^2(t) \\
x_2(t)
\end{bmatrix}\\
\mathbf{g}(\mathbf{y}(t)) &=& \mathbf{f}^{-1} (\mathbf{y}(t)) &=& \begin{bmatrix} y_1(t) - y_2^2(t) \\ y_2(t) \end{bmatrix} .
\end{aligned}
\end{equation}
Figure~\ref{fig:initial_data_nonlinear} shows the data from Figure~\ref{fig:initial_data} transformed by the function $\mathbf{f}$ in \eqref{eq:nonlinear_function} and colored by time.
It is important to note that this is a difficult class of problem in practice, as none of the observed
variables are purely fast or slow, and the observed system appears, at first inspection, to possess no separation
of time scales.
For this example, the analytical covariance and inverse covariance are
\begin{equation}
\begin{aligned}
\mathbf{C}(\mathbf{x}(t)) =&
\frac{1}{\epsilon}
 \begin{bmatrix}
\epsilon + 4x_2^2(t) & 2x_2(t) \\
2x_2(t) & 1
\end{bmatrix}\\
\mathbf{C}^{\dagger}(\mathbf{x}(t)) =&
\begin{bmatrix}
1 & -2 x_2(t) \\
-2 x_2(t) & \epsilon+ 4 x_2^2(t)
\end{bmatrix} .
\end{aligned}
\end{equation}

The fast and slow variables are now coupled through the function $\mathbf{f}$, and the Euclidean distance is not informative about the fast {\em or} the slow variables.
We need to use the Mahalanobis distance to obtain a parametrization that is consistent with the underlying fast-slow dynamics.

\subsubsection{Errors in Mahalanobis distance}

We can bound the Mahalanobis distance by the eigenvalues of $\mathbf{C}^{\dagger}$,
\begin{equation}
\lambda_{C^{\dagger},1} \| \mathbf{y}(t_2) - \mathbf{y}(t_1) \|_2^2
\le
\| \mathbf{y}(t_2) - \mathbf{y}(t_1) \|^2_M
\le
\lambda_{C^{\dagger},2} \| \mathbf{y}(t_2) - \mathbf{y}(t_1) \|_2^2
\end{equation}
where $\lambda_{C^{\dagger},1} \le \lambda_{C^{\dagger},2}$ are the two eigenvalues of $\mathbf{C}^{\dagger}$.
Therefore, for the example in \eqref{eq:nonlinear_function}, we have
\begin{equation}
E_M(\mathbf{y}(t_1), \mathbf{y}(t_2)) = - (y_2(t_2) - y_2(t_1))^4.
\end{equation}

\begin{figure}[t]
\centering
\begin{subfigure}{0.4\textwidth}
\centering
\epsfig{height=2in, file=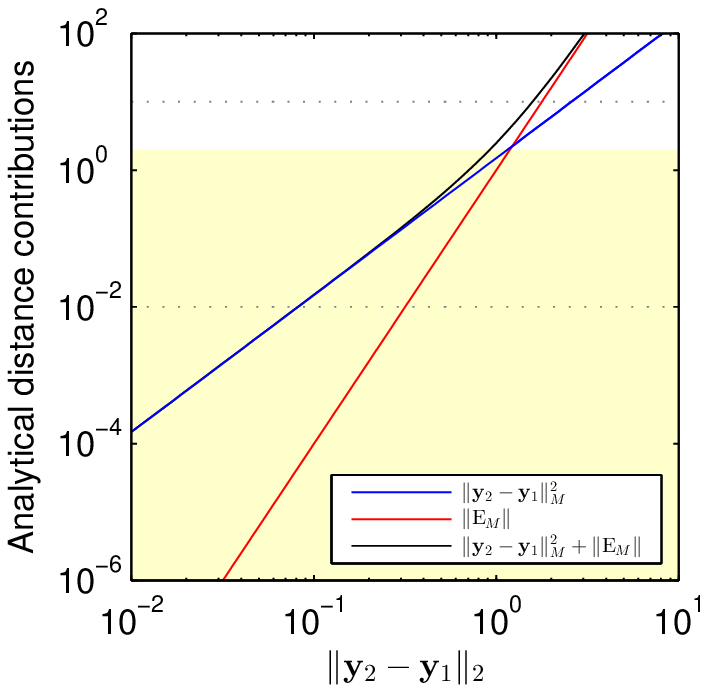}
\caption{}
\end{subfigure}
\begin{subfigure}{0.4\textwidth}
\centering
\epsfig{height=2in, file=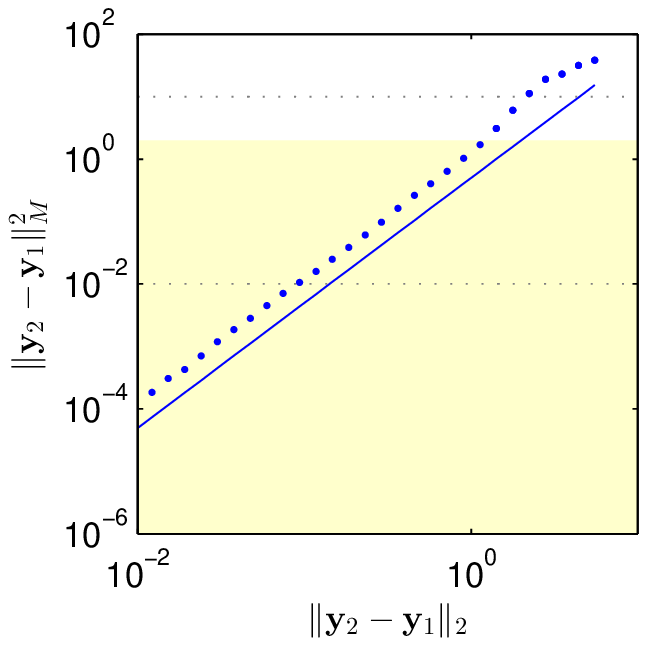}
\caption{}
\end{subfigure}

\begin{subfigure}{0.4\textwidth}
\centering
\epsfig{height=2in, file=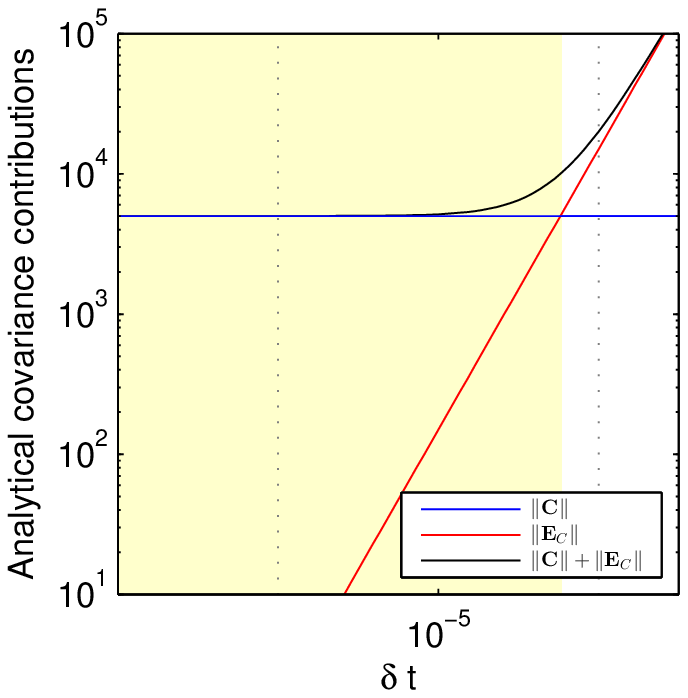}
\caption{}
\end{subfigure}
\begin{subfigure}{0.4\textwidth}
\centering
\epsfig{height=2in, file=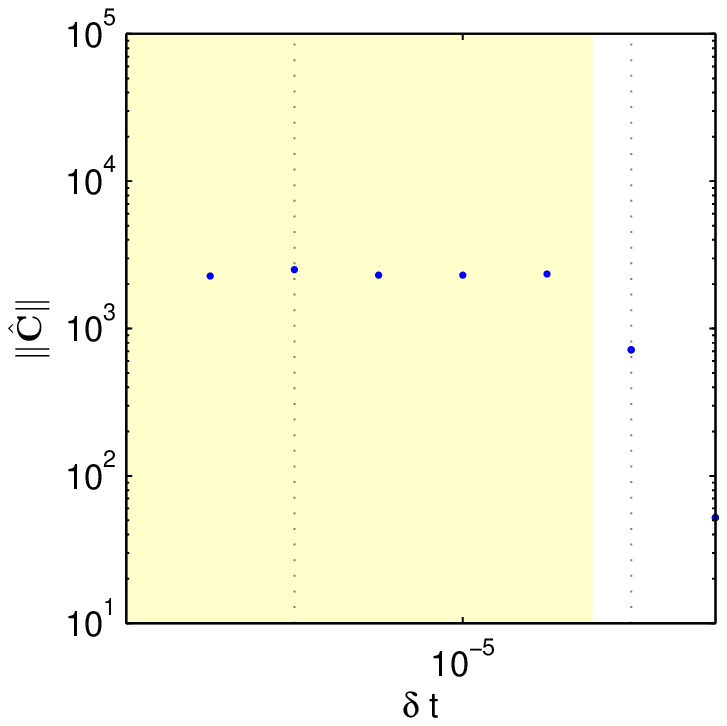}
\caption{}
\end{subfigure}
\caption{Errors in the Mahalanobis distance and the covariance estimation for the nonlinear example in \eqref{eq:nonlinear_function}. (a) The analytical expressions for the contributions to the distance approximation as a function of $\| \mathbf{y}_2 - \mathbf{y}_1\|_2$. (b) The average estimated Mahalanobis distance $\| \mathbf{y}_2 - \mathbf{y}_1\|^2_M$ as a function of the distance $\| \mathbf{y}_2 - \mathbf{y}_1\|_2$. The line $\| \mathbf{y}_2 - \mathbf{y}_1\|^2_M = \| \mathbf{y}_2 - \mathbf{y}_1\|_2$ is shown for reference. The yellow region indicates the range in which $\sigma_{kernel}$ should be chosen. (c) The analytical expressions for the contributions to the covariance as a function of $\delta t$. (d) The average estimated covariance $\| \hat{\mathbf{C}} \|$ as a function of $\delta t$. The yellow region indicates the range of $\delta t$ over which the errors in the estimated covariance are small relative to the norm of the covariance. }
\label{fig:cov_error_nonlinear}
\end{figure}

Figure~\ref{fig:cov_error_nonlinear}(a) shows $\| \mathbf{y}(t_2) - \mathbf{y}(t_1) \|^2_M$ and $| E_M |$ as a function of $\| \mathbf{y}(t_2) - \mathbf{y}(t_1) \|_2$.
The Mahalanobis distance is an accurate approximation to the true intrinsic distance $\| \mathbf{z}(t_2) - \mathbf{z}(t_1) \|_2$ when $|E_M| \ll \| \mathbf{y}(t_2) - \mathbf{y}(t_1) \|^2_M$ (the shaded yellow region in the plot indicates where $|E_M| < \| \mathbf{y}(t_2) - \mathbf{y}(t_1) \|^2_M$).
We want to choose $\sigma_{kernel}^2$ in a regime where $|E_M(\mathbf{y}(t_1), \mathbf{y}(t_2))| \ll \| \mathbf{y}(t_2) - \mathbf{y}(t_1) \|^2_M$, so that the distances we utilize in the diffusion maps calculation are accurate.
We can find such a regime empirically by plotting $\| \mathbf{y}(t_2) - \mathbf{y}(t_1) \|^2_M$ as a function of $\| \mathbf{y}(t_2) - \mathbf{y}(t_1) \|_2$, and assessing when the relationship deviates from quadratic.
This is shown in Figure~\ref{fig:cov_error_nonlinear}(b), and the deviation from quadratic behavior is consistent with the intersection of the analytical expressions plotted in Figure~\ref{fig:cov_error_nonlinear}(a).

Figures~\ref{fig:colored_data_nonlinear_cases}(a) and (b) show the data from Figure~\ref{fig:initial_data_nonlinear}, colored by $\phi_1$ for two different values of $\sigma_{kernel}$.
The corresponding values of $\sigma_{kernel}^2$ are indicated by the dashed lines.
When $\sigma_{kernel}^2$ corresponds to a region where $ |E_M | \ll \|\mathbf{y}(t_2) - \mathbf{y}(t_1) \|_M^2$, $\phi_1$ is well correlated with the slow variable.
However, when $\sigma_{kernel}^2$ corresponds to a region where $|E_M | \gg \|\mathbf{y}(t_2) - \mathbf{y}(t_1) \|_M^2$, the slow variable is no longer recovered.

\subsubsection{Errors in covariance estimation}

From \eqref{eq:cov_error}, we find that, for the example in \eqref{eq:nonlinear_function},
\begin{equation}
\small
\begin{aligned}
E_{C,11} (\mathbf{x}(t), \delta t)
=&
\frac{2 \delta t}{\epsilon^2}
- \frac{8 x_2(t)}{\epsilon^2 \delta t} \mathbb{E} \left[ \int_t^{t+\delta t} \left( \int_{s_2}^{t+\delta t} 2 x_2(s_1) ds_1
+  \int_t^{s_2} x_2(s_1) ds_1 \right) ds_2\right]  + \mathcal{O} (\delta t^{3/2}) \\
E_{C, 12} (\mathbf{x}(t), \delta t)
= &
E_{C, 21} (\mathbf{x}(t), \delta t)\\
=&
- \frac{x_2(t) \delta t}{\epsilon^2}
- \frac{2}{\epsilon^2 \delta t} \mathbb{E} \left[ \int_t^{t+\delta t} \left( \int_{s_2}^{t + \delta t} 2 x_2(s_1) ds_1 + \int_t^{s_2} x_2(s_1) ds_1 \right) ds_2 \right] + \mathcal{O} (\delta t^{3/2})\\
E_{C, 22} (\mathbf{x}(t), \delta t)
=&
-\frac{\delta t}{\epsilon^2} + \mathcal{O} (\delta t^{3/2})
\end{aligned}
\end{equation}
The error in the covariance is $\mathcal{O} \left( \frac{\delta t}{\epsilon^2} \right)$.
As expected, the error grows with increasing $\delta t$.
We can also see the explicit dependence of the covariance error on the time scale separation $\epsilon$; larger time scale separation results in a larger covariance error, as a more refined simulation burst is required to estimate the covariance of the fast directions.
$\|\mathbf{C} \|$ and $\| \mathbf{E}_C\|  $ are plotted as a function of $\delta t$ in Figure~\ref{fig:cov_error_nonlinear}(c); the shaded yellow portion denotes the region where $\| \mathbf{E}_C \| < \| \mathbf{C} \|$.
As in the previous example, we can empirically find where $\| \mathbf{E}_C \| \ll \| \mathbf{C} \|$  by plotting $\| \hat{\mathbf{C}} \|$ as a function of $\delta t $ and looking for a knee in the plot.
These results are shown in Figure~\ref{fig:cov_error_nonlinear}(d).

Figures~\ref{fig:colored_data_nonlinear_cases}(a) and (c) show the data from Figure~\ref{fig:initial_data_nonlinear}, colored by $\phi_1$ for two different values of $\delta t$.
The corresponding values of $\delta t$ are indicated by the dashed lines in Figure~\ref{fig:cov_error_nonlinear}(c)~and~(d).
When $\delta t$ corresponds to a region where $\|\mathbf{E}_C \| \ll \| \mathbf{C} \|$, the slow variable is recovered by the first diffusion maps coordinate.
However, when $\delta t$ corresponds to a region where $\|\mathbf{E}_C \| \gg \| \mathbf{C} \|$, the slow variable is no longer recovered.

\def \figheight {1.6in}

\begin{figure}[t]

\begin{subfigure}{\textwidth}
\epsfig{height=\figheight, file=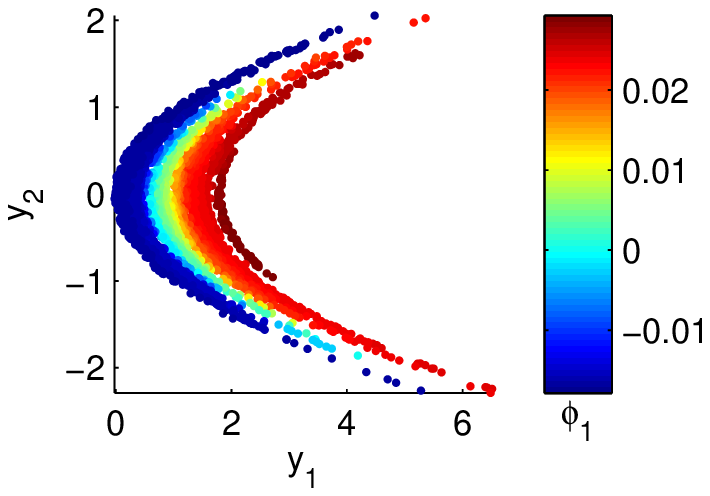}
\hfill
\epsfig{height=\figheight, file=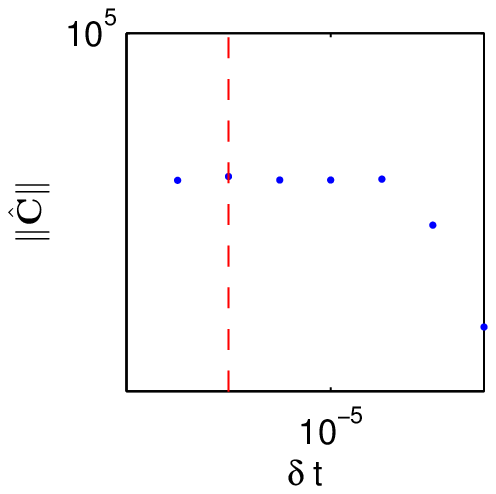}
\hfill
\epsfig{height=\figheight, file=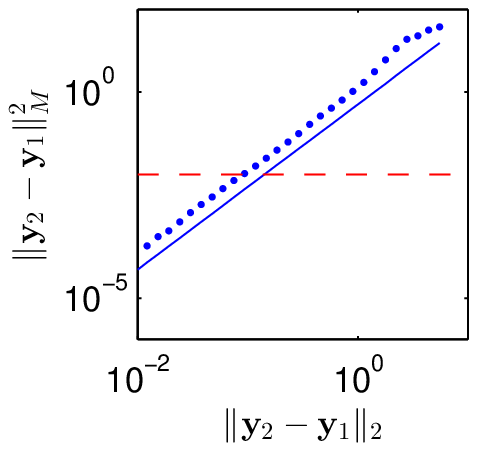}
\end{subfigure}

\begin{subfigure}{\textwidth}
\epsfig{height=\figheight, file=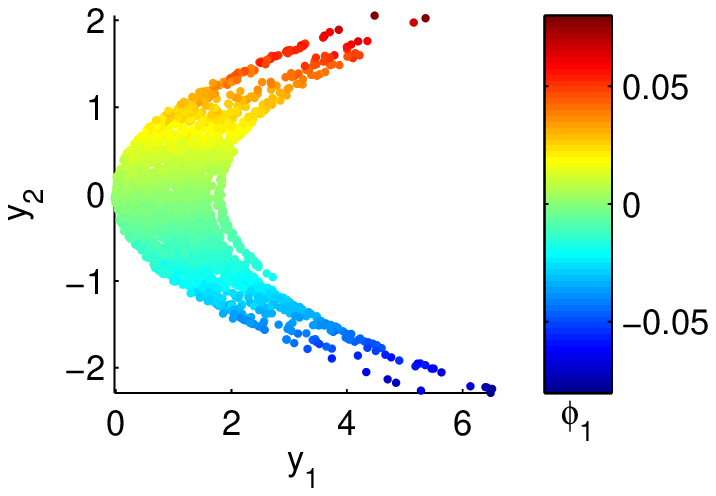}
\hfill
\epsfig{height=\figheight, file=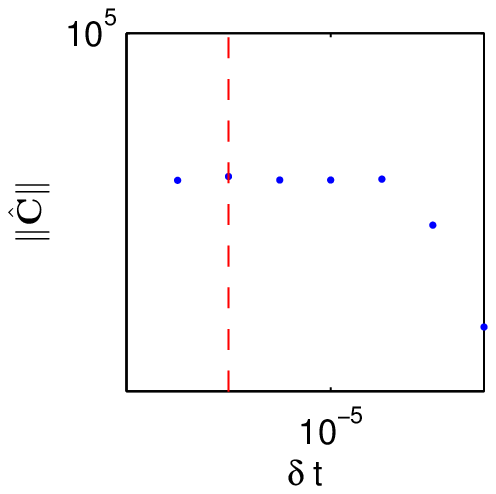}
\hfill
\epsfig{height=\figheight, file=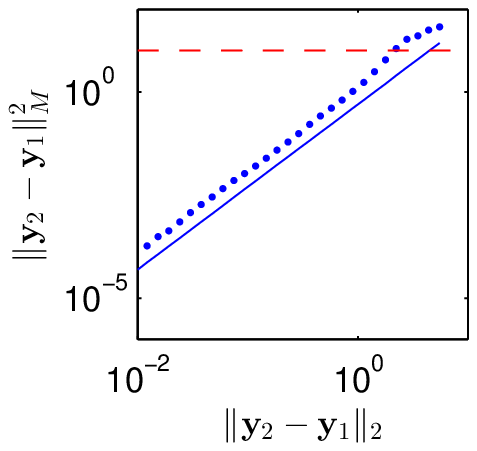}
\end{subfigure}

\begin{subfigure}{\textwidth}
\epsfig{height=\figheight, file=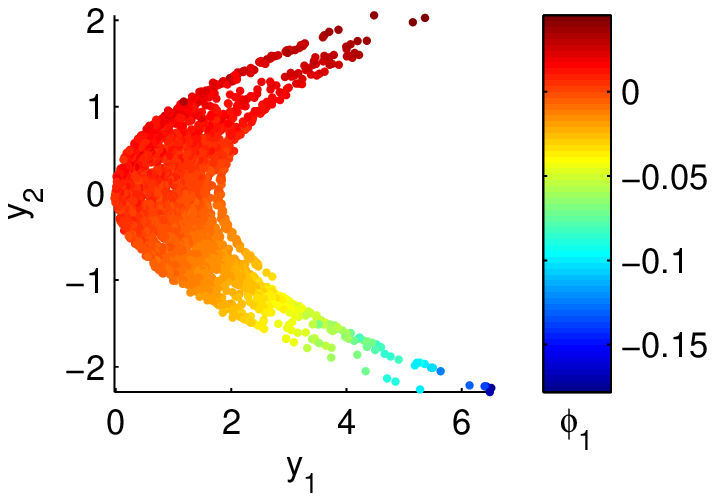}
\hfill
\epsfig{height=\figheight, file=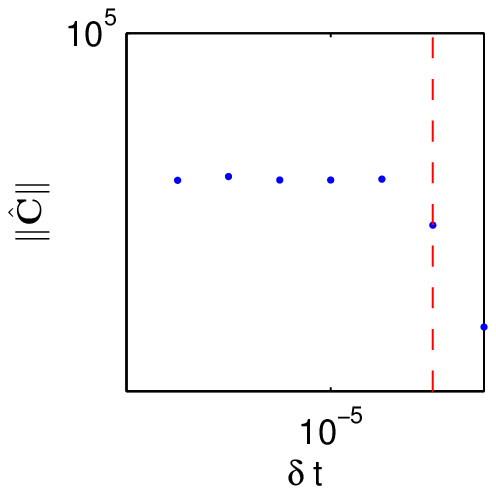}
\hfill
\epsfig{height=\figheight, file=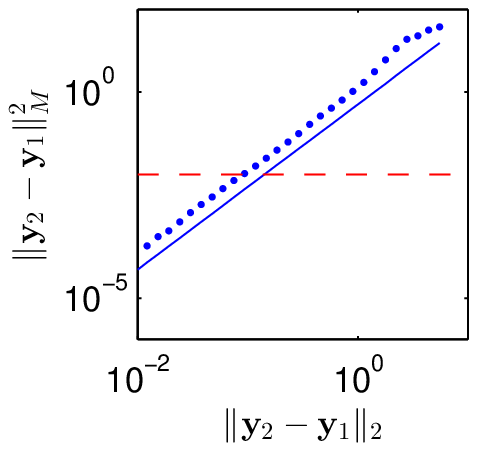}
\end{subfigure}

\caption{Data from Figure~\ref{fig:initial_data_nonlinear}, colored by the first diffusion maps variable $\phi_1$ using the Mahalanobis distance for three different parameter settings. The relevant values of $\delta t$ and $\sigma_{kernel}$ are indicated by the red dashed lines on the corresponding plots.  (Row~1) $\delta t = 10^{-7}$ and $\sigma_{kernel}^2 = 10^{-2}$. Note that the parametrization is one-to-one with the slow variable.  (Row~2) $\delta t = 10^{-7}$ and $\sigma_{kernel}^2 = 10^{1}$. We do not recover the slow variable because $\sigma_{kernel}$ is too large. (Row~3)  $\delta t = 10^{-3}$ and $\sigma_{kernel}^2 = 10^{-2}$. We do not recover the slow variable because $\delta t$ is too large.  }
\label{fig:colored_data_nonlinear_cases}
\end{figure}

\section{Conclusions}

We have presented a methodology to compute a parametrization of a data set which respects the slow variables in the underlying dynamical system.
The approach utilizes diffusion maps, a kernel-based manifold learning technique, with the Mahalanobis distance as the metric.
We showed that the Mahalanobis distance collapses the fast directions within a data set, allowing for successful recovery of the slow variables.
Furthermore, we showed how to estimate the covariances (required for the Mahalanobis distance) directly from data.
A key point in our approach is that the embedding coordinates we compute are not only insensitive to the fast variables, but are also invariant to nonlinear observation functions.
Therefore, the approach can be used for data fusion: data collected from the same system via different measurement functions can be combined and merged into a single coordinate system.

In the examples presented, the initial data came from a single trajectory of a dynamical system, and the local covariance at each point in the trajectory was estimated using brief simulation bursts.
However, the initial data need not be collected from a single trajectory, and other sampling schemes could be employed.
Brief time series are required to estimate the local covariances, but given a simulator, one could reinitialize brief simulation bursts which are sufficiently short and refined from each sample point.

In our examples, we controlled the time scale of sampling and could therefore set the time scale over which to estimate the covariance and the simulation time step arbitrarily small.
However, in some settings, such as previously collected historical data, it is not uncommon to have a fixed sampling rate and be unable to reinitialize simulations.
In such cases, it is possible that we cannot find an appropriate kernel scale given the fixed $\delta t$ such that we can accurately recover the slow variables.
For these cases, the data cannot be processed as given,
and it is necessary to construct intermediate observers,
such as histograms, Fourier coefficients, or scattering transform coefficients \cite{mallat2012group, talmon2014intrinsic, talmon2014manifold}.
Such intermediates are more complex statistical functions than simple averages and can capture additional structure within the data.
They also reduce the effects of noise and permit a larger time step.
However, constructing such intermediates often requires additional {\em a priori} knowledge about the system dynamics and noise structure.

Clearly, in our analysis, we have ignored the finite sampling effects in our estimation.
In reality, both the number of samples used to estimate the covariances, as well as the density of sampled points on the manifold, affect the recovered parametrization and provide additional constraints on $\delta t$ and $\sigma_{kernel}$.
Future work involves extending our analysis to the finite sample case, and providing guidelines for the amount of data required to apply our methodology.

The methods presented here provide a bridge between traditional data mining and multiple time scale dynamical systems.
With this interface established, one can now consider using such data-driven methodologies to extract reduced models (either explicitly, or implicitly via an equation-free framework \cite{erban2006gene, kevrekidis2004equation, kevrekidis2003equation,  kevrekidis2009equation}) which also respect the underlying slow dynamics and geometry of the data.
Such reduced models hold the promise of accelerated analysis and reduced simulation of dynamical systems whose effective dynamics are obscure upon simple inspection.

\bibliographystyle{siam}
\bibliography{../../../../references/references}

\end{document}